\documentclass[twocolumn,final]{autart}    
\usepackage{enumerate}
\usepackage{apalike}
\usepackage{amssymb}
\usepackage{amsmath}
\usepackage{graphicx}          

\newcommand{\ba}{\begin{array}}
\newcommand{\ea}{\end{array}}
\newcommand{\noi}{\noindent}
\newcommand{\Frac}{\displaystyle\frac}
\newcommand{\tp}{^\mathsf{T}}
\newcommand{\V}{\mathcal{V}}
\newcommand{\I}{\boldsymbol{\mathcal{I}}}
\newcommand{\B}{\boldsymbol{\mathcal{B}}}
\newcommand{\nn}{\nonumber}
\newcommand{\eps}{\epsilon}

\begin{document}

\begin{frontmatter}

\title{A New Position Control Strategy for VTOL UAVs using IMU and GPS measurements} 

\thanks[footnoteinfo]{This paper was not presented at any IFAC
meeting. Corresponding author A. Tayebi. Tel. +1 (807) 343-8597. Fax +1 (807) 766-7243.}

\author[a]{Andrew Roberts}\ead{arober88@uwo.ca},    
\author[a,b]{Abdelhamid Tayebi}\ead{atayebi@lakeheadu.ca}               

\address[a]{Department of Electrical and Computer Engineering, University of Western Ontario, London, Ontario, Canada, N6A 3K7.}             
\address[b]{Department of Electrical Engineering, Lakehead University, Thunder Bay, Ontario, Canada P7B 5E1}

\begin{keyword}                           
VTOL, UAV, vector measurments, position control,             
\end{keyword}                             

\begin{abstract}                          
We propose a new position control strategy for VTOL-UAVs using IMU and GPS measurements. Since there is no sensor that measures the attitude, our approach does not rely on the knowledge (or reconstruction) of the system orientation as usually done in the existing literature. Instead, IMU and GPS measurements are directly incorporated in the control law. An important feature of the proposed strategy, is that the accelerometer is used to measure the apparent acceleration of the vehicle, as opposed to only measuring the gravity vector, which would otherwise lead to unexpected performance when the vehicle is accelerating (i.e. not in a hover configuration). Simulation results are provided to demonstrate the performance of the proposed position control strategy in the presence of noise and disturbances.
\end{abstract}

\end{frontmatter}

\section{Introduction}

The design of position controllers for Vertical take-off and landing (VTOL) unmanned airborne vehicles (UAVs) has been the focus of several research groups, which has resulted in significant breakthroughs in this field, for example see \cite{Abdelkader2010}, \cite{Aguiar2007}, \cite{Frazzoli2000}, \cite{Hauser1992}, \cite{Hua2009}, \cite{Pflimlin2007}  and \cite{Roberts2011}. Existing position controllers, usually require that the system states are accurately known or measured, namely the position, linear velocity, angular velocity and the orientation. For outdoor applications a global positioning system (GPS) mounted to the system can be used to provide the position and velocity measurements, while the angular velocity is obtained using a gyroscope which is included in the inertial measurement unit (IMU) in addition to an accelerometer and a magnetometer. However, there does not exist any sensor that provides directly the orientation of a rigid body. Motivated by this problem, the study of \emph{rigid-body attitude estimation} has seen substantial breakthroughs due to the efforts of the research community (see, for instance, \cite{Mahony2008} ). However, we are not aware of any work in the literature, providing a rigorous results for the combination of an attitude observer and a position controller for VTOL-UAVs.\\
To address this shortcoming, there has been some effort to design position control algorithms which do not directly require the measurement of the system attitude. For example, in \cite{Roberts_IFAC_2011} the authors propose a position control law which utilizes a number of \textit{vector measurements} as a means to eliminate the requirement for the attitude measurement. By vector measurements we are referring to the body-referenced measurements of vectors whose coordinates are known in the inertial frame. Since the vector measurements contain information about the system orientation, it has been shown that they can be applied directly to the position controller thereby eliminating the need of the observer completely. Consequently, the resulting vector-measurement-based position control laws do not require the direct measurement of the system attitude, nor do they require an attitude observer which provides practitioners with a simpler, reduced order closed loop system, with accompanying proofs for stability.\\
Unfortunately, the vector-measurements based position control strategy can be susceptible to a problem associated with the lack of sensors which can provide suitable vector measurements. This shortcoming stems from the fact that the two sensors most commonly used to provide vector measurements are the magnetometer and accelerometer, used to provide body-referenced measurements of the Earths magnetic field and gravity vector, respectively. However, in order to satisfy the requirement that the accelerometer provides a measurement of the gravity vector only, one must assume that the body-fixed frame is non-accelerating. It is clear that this condition is not guaranteed to be satisfied in some applications involving VTOL UAVs. \\
Of course, this practical limitation is relevant to both vector-measurement-based position controllers and attitude observers which use accelerometers. Fortunately, this limitation has led to the development of a new class of attitude observers which uses the accelerometer (and magnetometer) to provide vector-measurements. This type of observer acknowledges the fact that the accelerometer measures a combination of the gravity vector and the acceleration of the rigid-body in the body-fixed frame. This combination of the gravity vector and linear acceleration in the inertial reference frame is commonly referred to as the \emph{apparent acceleration}. This inertial vector violates the requirement of many of the vector-measurement-based attitude observers, since the system acceleration is not known in the inertial frame of reference. In order to deal with the fact that the inertial vector is unknown, this type of attitude observer uses the velocity of the rigid-body (assumed to be measurable using, for instance, a GPS) in addition to the signals obtained from an IMU. These attitude observers, which are often referred to as \emph{velocity-aided attitude observers}, can be found in \cite{Bonnabel2008}, \cite{Martin2008} and \cite{Martin2010} with local stability proofs , and in \cite{Hua2010} with almost semiglobal stability results. \\
In this paper we propose a new position control approach which obviates the requirement of the system attitude measurement by using the vector measurements directly in the control law. We specifically use a magnetometer and accelerometer
to provide the two vector measurements. The accelerometer is used to measure the system \emph{apparent acceleration}, rather than the gravity vector only.
Using our proposed approach, we show that, upon a suitable choice of the control gains, all system states remain bounded,
and the system position converges to a constant reference position. Our proposed control strategy 1)  does not require direct measurement of the system attitude; 2) does not require the use of an attitude observer; 3) uses an accelerometer to provide a vector measurement without limiting the motion of the system to a near-hover state.
%


\section{Background}
In this section we present some of the necessary mathematical details we use throughout the paper. In section \ref{section:attitude_representation} we describe two commonly used attitude representations (rotation matrices and unit-quaternion). In section \ref{section:math_background} we define functions which are necessary in developing the proposed control laws.

\subsection{Attitude Representation}
\label{section:attitude_representation}
To represent the orientation of the aircraft (rigid-body), we define two reference frames:  An inertial frame $\I$, which is rigidly attached the Earth (assumed flat), and a body frame $\B$ which is rigidly attached to the aircraft center of gravity (COG). The orthonormal basis of $\mathcal{B}$ is taken such that the $x$ axis is directed towards the front of the aircraft (or rigid body), the $y$ axis is taken towards the starboard (right) side, and the $z$ axis is directed downwards (opposite the direction of the system thrust).\\
Throughout the paper we often refer to the \emph{orientation} of the rigid-body, by which we mean the relative angular position of $\B$ with respect to $\I$. The goal of the attitude representation is to mathematically describe the orientation of the rigid-body.
The unit-quaternion, which is a unit vector on $\mathbb{R}^4$, is given by $Q = \left(\eta,q\right) \in \mathbb{Q}$, where $\eta \in \mathbb{R}$ is the quaternion-scalar
and $q \in \mathbb{R}^3$ is the quaternion-vector, and $\mathbb{Q}$ is the set of unit-quaternion defined by
\begin{equation}
\mathbb{Q}\equiv \left\{ Q \in \mathbb{R}\times \mathbb{R}^3 ,~ \| Q \| = 1 \right\}.
\end{equation}
Let $Q_1 = (\eta_1,q_1) \in \mathbb{Q}$, $Q_2 = (\eta_2,q_2) \in \mathbb{Q}$ denote two unit-quaternion; then the \emph{quaternion product} of $Q_1$ and $Q_2$, denoted by $Q_3 = (\eta_3,q_3) \in \mathbb{Q}$ is defined by the following operation
\begin{equation}
Q_3 = Q_1 \odot Q_2 =
\left(
\ba{cc}
\eta_1 \eta_2 - q_1\tp q_2,
&
\eta_1 q_2 + \eta_2 q_1 + S(q_1)q_2
\ea
\right).
\end{equation}
The set of unit-quaternion $\mathbb{Q}$ forms a group with the quaternion multiplication operation $\odot$,
with the quaternion inverse $Q^{-1} = (\eta,-q)$, and identity element $(1,\mathbf{0})=Q^{-1}\odot Q = Q \odot Q^{-1} $. \
The unit-quaternion is an over-parameterization of the the special group of orthogonal matrices of dimension three $SO(3)$, defined as
\begin{equation}
SO(3) \equiv \left\{ R \in \mathbb{R}^{3 \times 3}, |R| = 1, R\tp R = R R\tp = I \right\},
\end{equation}
that is, the transformation from the quaternion space $\mathbb{Q}$ to $SO(3)$, given by the following \emph{Rodrigues} formula:
\begin{equation}
R(Q) =  I + 2 S(q)^2 - 2 \eta S(q).
\label{q2r}
\end{equation}
where $S(\cdot)$ is a skew-symmetric matrix given in the next section, is a two-to-one map, \textit{i.e.,} $R(Q)=R(-Q)$.

\subsection{Skew symmetric matrices and bounded functions}
\label{section:math_background}
Let $x,y \in \mathbb{R}^3$. We define the skew-symmetric matrix $S(x)$  such that $S(x)y=x \times y$, where $\times$ denotes the vector cross product.
%
%
Several useful properties of this skew-symmetric matrix are given below:
\begin{align}
S(x)^2  &=  x x\tp - x\tp x I_{3\times 3},
\label{property:skew:squared}
\\
S(R x ) &=  R S(x) R\tp, \quad R \in SO(3),
\label{property:skew:so3}
\\
S(x)y & =  - S(y)x ~ = ~x \times y,
\label{property:skew:negate}
\\
\lambda \left( S(x)^2 \right)
&=
\left[\ba{ccc}
0, & -\|x\|^2, & -\|x\|^2
\ea\right],
\label{property:skew:eigenvalues}
\end{align}
where $\lambda(M)$ denotes the eigenvalues of the matrix $M \in \mathbb{R}^{3\times 3}$.\\
Consider the bounded, differentiable function, denoted as $h(\cdot) : \mathbb{R}^3 \rightarrow \mathbb{R}^3$, which satisfies the following properties:
\begin{equation}
\ba{cc}
u \tp h(u)  > 0 & \forall u \in \mathbb{R}^3, \|u\| \in (0,\infty),
\\
\left.\ba{c}
0 \leq \| h(u) \| < 1
\\
0 < \| \phi(u) \| \leq 1
\ea\right\}
&
\forall u \in \mathbb{R}^3, \|u\| \in [0,\infty),
\ea
\end{equation}
where  $\phi(u) := \frac{\partial}{\partial u} h(u)$. Throughout the paper we make use of one particular example given by $h(u) = \left( 1 + u\tp u \right)^{-1/2}u $. Using this definition for $h(\cdot)$ one can derive the expression
$
\phi(u)
=
 (1 + u\tp u )^{-3/2} ( I_{3\times 3} - S(u)^2 ).
$.
\section{Position Control Using GPS and IMU Measurements}
Using the above mathematical background, we will now proceed to formulate the problem and define the position control laws. A number of steps are taken which are grouped into various sections. In section \ref{section:equations_motion} we define the system model. In section \ref{section:problem} we formulate the problem and state some necessary assumptions. Section \ref{section:attitude_extract} provides an attitude extraction algorithm which allows us to specify a desired system attitude based upon the values of the position and velocity error, and section \ref{section:attitude_error} defines attitude error functions. Finally, in section \ref{section:controller} we describe the position control laws.
\subsection{Equations of Motion}
\label{section:equations_motion}
To model the system translational dynamics, we let $p, v \in \mathbb{R}^3$ denote the position and velocity, respectively, of the vehicle COG expressed in the inertial frame $\mathcal{I}$. For this problem we assume that the body-referenced angular velocity vector $\omega$ is available as a control input. We consider the following VTOL UAV model:
\begin{align}
 \dot{p}
 &= v,
 \label{control:p_dot}
 \\%
 \dot{v}
 &= \mu  + \delta, \quad \mu = g e_3 - u_t R\tp e_3,
\label{control:v_dot}
\\
\dot Q
&=  \frac{1}{2} \left[\ba{c}
-q \tp \\
\eta I_{3\times 3} + S(q)
\ea \right] \omega,
\label{control:Q_dot}
\end{align}

\noi where $u_t = T/m_b$, $T$ is the system thrust, $m_b$ is the system mass, $e_3 = \mbox{col}\left[0,0,1\right]$, $g$ is the gravitational acceleration, and $\delta$ is a disturbance which is dependent on aerodynamic drag forces. The control input of the system is defined as $u = \left[ u_t, \omega \right]\tp$. The system output is defined as $y = \left[ p,v, b_1, b_2\right]\tp$ where $b_2$ is the signal obtained using an accelerometer, $b_1 = R r_1$ is a signal obtained using a magnetometer, and $r_1$ is the magnetic field of the surrounding environment (assumed constant). Note that the system attitude $R$ (or $Q$) is not assumed to be a known output of the system.

We consider a well-known model for the accelerometer model (which includes forces due to linear acceleration $\dot v$) which is given by
\begin{equation}
\ba{l}
b_2 =
R \left( \dot v - g e_3 \right)
 =
- u_t e_3 + R \delta
=  R r_2,
\label{control:accel_model}
\ea
\end{equation}
where $r_2$ is the inertial referenced system apparent acceleration, which satisfies
\begin{equation}
\dot v
~=~
 g e_3 + r_2,
\qquad
r_2
~= ~
 - u_t R\tp e_3 + \delta.
\label{control:v_dot:a}
\end{equation}

In the development of attitude observers, it is often assumed that the system is near hover (or $\dot v \approx 0$) in order to assume that the accelerometer measures the direction of the gravity vector. Also, in most situations the aerodynamic disturbance vector $\delta$ is not included in the model. However, for the VTOL UAV model, one can easily see that if the aerodynamic disturbance is neglected, or we assume that $\delta \approx 0$, then the accelerometer signal provides the measurement $b_2 = - u_t e_3$, which is the constant vector $e_3$ multiplied by the system thrust. In this case the use of the accelerometer seems trivial since its measurement is known a priori and does not contain any information about the system attitude. Therefore, we see that for the VTOL UAV system, the assumption that the accelerometer measures only the gravity vector may be a dangerous assumption which may lead to unexpected performance, even in the case where $\dot v \approx 0$. In fact, it seems that the utility of the accelerometer measurements is related to the measurement of the vector $\delta $ since the accelerometer measures $b_2 = - u_t e_3 + R \delta$. For this reason we believe that it is important to include a model of the aerodynamic disturbances.
\subsection{Problem Formulation}
\label{section:problem}
Let $p_r$ denote a desired reference position, which is assumed to be constant (or slowly-varying), and let $e_p = p - p_r$. Our main objective is to develop a control law for the system inputs $u_t$ and $\omega$, using the available system outputs $y = \left[ p,v,b_1,b_2\right]$, such that the system states $e_p$ and $v$ are bounded and $\lim_{t\to\infty}  e_p(t) = \lim_{t\to\infty} v(t)  = 0$. For the position control design we first require the following assumptions are satisfied.

\begin{assum}
There exist positive constants $c_1$ and $c_2$ such that
$
\| r_2(t) \| \leq c_1
$
and
$
\| \dot r_2(t) \| \leq c_2.
$
\label{assumption:r_2_constants}
\end{assum}
\begin{assum}
Given two positive constants, $\gamma_1$ and $\gamma_2$, there exists a positive constant $c_w(\gamma_1,\gamma_2)$ such that $c_w < \lambda_{\min}(W)$ where
$
W = -\gamma_1 S(r_1)^2 - \gamma_2 S(r_2)^2.
$
\label{assumption:lambda_min_w}
\end{assum}
The second assumption is satisfied if $r_2$ is non-vanishing and is not collinear to the magnetic field vector $r_1$. In the case where $r_2 = 0$, the system velocity dynamics become $\dot v = g e_3$ (which corresponds to the rigid body being in a free-fall state) which is not likely in normal circumstances. When this assumption is satisfied, it follows that $W$ is positive definite. Furthermore, if this assumption is satisfied, the value of $c_w>0$ can be arbitrarily increased by increasing the values of $\gamma_1$ and $\gamma_2$.\\
In addition to this assumption, we also require some conditions on the aerodynamic force vector $\delta$.
\begin{assum}[Aerodynamic Forces]
\label{assumption:aerodynamic_forces}
In light of the fact that the disturbance force $\delta$ is due to aerodynamic forces exerted on the vehicle we make the following simplifying assumptions:
\begin{enumerate}[(a)]
\item The aerodynamic disturbance $\delta$ is dissipative with respect to the system translational kinetic energy and satisfies $\delta \tp v \leq 0$.
\label{assum:aero:a}
\item The aerodynamic disturbance force $\delta$ is only dependent on the system translational velocity, and there exist a positive constant $c_1$ such that $\| \delta \| \leq c_1 \| v\|^2$
\label{assum:aero:b}
\item There exists positive constants $c_2$ and $c_3$ such that $\| \dot \delta \| < c_2 + c_3 \| v \|^3$.
    \label{assum:aero:c}
\end{enumerate}
\end{assum}
\noi Assumption \ref{assumption:aerodynamic_forces}(\ref{assum:aero:a}) and \ref{assumption:aerodynamic_forces}(\ref{assum:aero:b}) can be realized when the system is operating in an environment where the exogenous airflow is negligible (no wind). Assumption \ref{assumption:aerodynamic_forces}(\ref{assum:aero:c}) can be satisfied when the system geometry is sufficiently symmetrical such that the system aerodynamic forces do not significantly depend on the system orientation. Although this assumption may be reasonable for certain VTOL type aircraft, for example the ducted-fan, this assumption may not be the case with certain systems, for example fixed wing aircraft, where the system aerodynamics depend largely on the orientation of the vehicle.
Now that we have established the required assumptions, let us consider the model for the system acceleration from \eqref{control:v_dot}.
Due to the underactuated nature of this system, the translational acceleration is driven by the system thrust and orientation $\mu(u_t,R)$. That is, if $\mu$ was a control input, setting $\mu = -k_p e_p - k_v v$ would satisfy the objectives (since $v\tp\delta \leq 0$). However, since $\mu$ is a function of the system state, we define $\mu_d \in \mathbb{R}^3$ as the \emph{desired acceleration}, and introduce the new error signal
\begin{equation}
\tilde\mu = \mu - \mu_d.
\label{control:e_mu}
\end{equation}
\noi Subsequently, a new objective is to force $\tilde\mu \rightarrow 0$ in order to obtain the desired translational dynamics. Since the signal $\mu$ is dependent on the system thrust and attitude, based upon the value of the desired acceleration $\mu_d$ we wish to obtain a suitable \emph{desired attitude}, denoted as $Q_d = \left(\eta_d,q_d\right)\in \mathbb{Q}$, and system thrust $u_t$, such that the following equation is satisfied
\begin{equation}
\mu_d = g e_3 - u_t R_d\tp e_3,
\label{control:mu_d_defn}
\end{equation}
where $R_d = R\left(Q_d\right)$ is the rotation matrix corresponding to the unit-quaternion $Q_d$, as defined by (\ref{q2r}). An extraction method satisfying these requirements, which has been previously given in \cite{Roberts2011}, is described in the following section.
\subsection{Desired Attitude and Thrust Extraction}
\label{section:attitude_extract}
In this section, given a value of the desired acceleration $\mu_d$, we seek to obtain the value of the desired orientation $R_d$ (or equivalently in terms of the unit-quaternion $Q_d$) such that equation (\ref{control:mu_d_defn}) is satisfied. To solve this problem we use an attitude and thrust extraction algorithm which has been previously proposed by \cite{Roberts2011}: Given $\mu_d$ where $\mu_d \notin L$,
\begin{equation}
L = \{ \mu_d \in \mathbb{R}^3; \mu_d = \mbox{col}[0,0,\mu_{d3}]; \mu_{d3} \in [g,\infty)\},
\label{L_singularity}
\end{equation}
\noi then, one solution for the thrust $u_t$ and attitude $Q_d = \left(\eta_d, q_d\right)$ where $R_d = R(Q_d)$, which satisfies (\ref{control:mu_d_defn}) is given by
\begin{align}
u_t
&=  \| \mu_d - g e_3 \|,
\label{extract:thrust}
\\
\eta_d
&=  \left( \frac{1}{2} \left( 1 + \frac{g - e_3\tp \mu_d}{\|\mu_d - g e_3 \| } \right)\right)^{1/2},
\label{extract:eta_d}
\\
q_d
&=  \frac{1}{2 \| \mu_d - g e_3 \| \eta_d } S(\mu_d) e_3.
\label{extract:q_d}
\end{align}
The extracted attitude $Q_d$ has the time-derivative
\begin{equation}
\dot Q_d =  \frac{1}{2}
\left[\ba{c}
-q_d \tp \\
\eta_d I_{3\times 3} + S(q_d)
\ea \right] \omega_d,
\label{extract:Q_d_dot}
\end{equation}
where the \textit{desired angular velocity} $\omega_d$ is given by
\begin{eqnarray}
\omega_d &= & M(\mu_d) \dot \mu_d,
\label{extract:omega_d_extract}
\\
M(\mu_d) &= & \frac{1}{4 \eta_d^2 u_t^4} \left( - 4 S(\mu_d) e_3 e_3\tp
+ 4 \eta_d^2 u_t S(e_3) + 2 S(\mu_d)
 \right. \nn \\ && {} \left.
 - 2 e_3 \tp \mu_d S(e_3)  \right) S\left(\mu_d - g e_3 \right)^2.
 \label{extract:M}
 \end{eqnarray}
\subsection{Attitude Error}
\label{section:attitude_error}
To represent the relative orientation of the desired attitude $Q_d$ with respect to the actual attitude $Q$, we let $\tilde Q = (\tilde\eta, \tilde q) \in \mathbb{Q}$ and $\tilde R = R ( \tilde Q) \in SO(3)$ denote the unknown \emph{attitude error} which is defined by
\begin{equation}
\tilde Q = Q \odot Q_d^{-1},
\qquad
\tilde R = R(\tilde Q ) = R_d\tp R,
\label{attitude_error}
\end{equation}
\noi where $Q_d$ is the unit quaternion obtained using (\ref{extract:eta_d}) and (\ref{extract:q_d}). In light of $\dot Q$ and $\dot Q_d$, as defined by (\ref{control:Q_dot}) and (\ref{extract:Q_d_dot}), respectively, the time derivative of the attitude error is found to be
\begin{equation}
\dot{\tilde Q} = \frac{1}{2} \left[
\ba{c}
-\tilde q\tp
\\
\tilde\eta I + S(\tilde q)
\ea\right] \tilde\omega,
\quad
\dot{\tilde R} = - S(\tilde\omega) \tilde R,
\label{control:Q_tilde_dot}
\end{equation}
\begin{equation}
\tilde\omega = R_d\tp \left(\omega - \omega_d\right),
\label{control:omega_tilde}
\end{equation}
where $\omega_d$ is the \emph{desired angular velocity} as defined by (\ref{extract:omega_d_extract}). One of the objectives of the control design is to force the system orientation to the desired attitude, or in terms of the rotation matrices, to force $R \rightarrow R_d$ (and therefore $\mu \rightarrow \mu_d$), in order to obtain the desired translational dynamics. As mentioned in section \ref{section:attitude_representation}, this corresponds to two possible solutions for the unit-quaternion which are given by $\tilde Q = \left( \pm 1, \mathbf{0} \right)$. The multiplicity of equilibrium solutions is manageable since our objectives are satisfied for both values of the unit-quaternion.
\subsection{Position Controller}
\label{section:controller}
The position controller design is based upon a value of the desired system translational acceleration, which is specified by the virtual control law $\mu_d$. Using the calculated values for the position error $e_p = p - p_r$,  and the system velocity $v$, the value of the desired acceleration is obtained. This desired acceleration is directly related to a corresponding desired rigid-body orientation and thrust, denoted by $Q_d$ and $u_t$, respectively, which is obtained using the attitude and thrust extraction method described in section \ref{section:attitude_extract}. The desired attitude given in the $SO(3)$ parametrization, denoted as $R_d$, is subsequently obtained using $Q_d$ with (\ref{q2r}).
Since the system attitude is not known, we incorporate the use of a special filter which is driven by the value of the linear velocity $v$. We let $\hat v \in \mathbb{R}^3$ denote the filter state variable which corresponds to the system velocity $v$, and define the error function $\tilde v = v - \hat v$.\\
Although the system linear velocity is known, the use of the signal $\hat v$ through the error function $\tilde v$, for an appropriate choice of the estimation law $\dot{\hat v}$ can can be viewed as a function of the system acceleration in terms of the unknown signal $r_2$. Since this vector is known in the body fixed frame (measured using an accelerometer, $b_2 = R r_2$), the filter variable $\hat v$ through the error function $\tilde v$ can be used with the accelerometer to provide information related to the system attitude. After these steps, the remaining control design is focused on forcing the actual system attitude to the desired attitude using the control input $\omega$.
The proposed control law is given as follows:
\begin{align}
\omega
&=
M(\mu_d) \left( f_{\mu_d} - k_v \phi(v) R_d\tp \left( b_2 + u_t e_3 \right) \right) + \psi,
\label{design:omega}
\\
f_{\mu_d}
&=
- k_p \phi(e_p) v + k_v \phi(v) \left( k_p h(e_p) + k_v h(v) \right),
\\
\psi
&=
\gamma_1 S(R_d r_1) b_1 + \gamma_2 k_1 S\left(R_d \left( v - \hat v \right)\right)b_2,
\label{design:psi}
\\
\dot{\hat v}
&=
g e_3 + R_d\tp b_2 + k_1 \left( v - \hat v\right)
+ \frac{1}{k_1} R_d\tp S(b_2) \psi,
\label{design:v_hat_dot}
\\
\mu_d
&=
- k_p h(e_p) - k_v h(v),
\label{design:mu_d}
\end{align}
where $k_1,\gamma_1,\gamma_2>0$, $M(\mu_d)$ is the function defined by (\ref{extract:M}), $\phi(\cdot)$ is the bounded function defined in section \ref{section:math_background},  $u_t = \| \mu_d - g e_3\|$, $R_d = R(Q_d)$ and $Q_d = (\eta_d,~q_d)$ is obtained from the value of $\mu_d$ using the attitude extraction algorithm defined in section \ref{section:attitude_extract}.
\begin{thm}
Consider the system given by (\ref{control:p_dot})-(\ref{control:Q_dot}), where we apply the control laws $u_t = \| \mu_d - g e_3\|$ and $\omega$ as defined by (\ref{design:omega}), where $k_p>0$ and $k_v>0$ are chosen such that $k_p + k_v < g$. Let assumptions \ref{assumption:lambda_min_w} and \ref{assumption:aerodynamic_forces} be satisfied. Then the system thrust $u_t$ is bounded and non-vanishing such that
 \begin{equation}
    0
    <
    \underline c_t
    \leq
    u_t(t)
    \leq
     \bar c_t,
    \quad
    \underline c_t
    =
    g - k_p - k_v,
    \quad
    \bar c_t
    =
    g + k_p + k_v,
    \label{proof:thrust_bound}
\end{equation}
and for all initial conditions $\tilde\eta(t_0) \neq 0$ (or equivalently $\|\tilde q(t_0)\| \neq 1$), there exists positive constants $\bar\gamma_1,\bar\gamma_2,\kappa_1>0$ such that for $\gamma_1 > \bar\gamma_1$, $\gamma_2 > \bar \gamma_2$, $k_1 > \kappa_1$, the system states $e_p$ and $v$ are bounded and $\lim_{t\to\infty} e_p(t) = \lim_{t\to\infty} v(t)  = \mathbf{0}$.
\end{thm}
%
%
%
%
\textbf{Proof:}
We begin by first proving the upper and lower bounds on the thrust control input $u_t$. Since the function $h(\cdot)$ is bounded by unity, the norm of the virtual control law $\mu_d$ is bounded by $\|\mu_d\| < k_p + k_v$. Since the thrust control input is given by $u_t = \| \mu_d - g e_3\|$, and $k_p$ and $k_v$ are chosen such that $k_p + k_v < g$, one easily arrives at the lower and upper bounds for $u_t$ described in the theorem. A nice consequence of the boundedness of $u_t$, is that the function $M(\mu_d)$ defined by (\ref{extract:M}), which is used in the expression for the desired angular velocity $\omega_d$, is also bounded. In fact, in \cite{Roberts2011} the authors show that the norm of this matrix satisfies
\begin{equation}
\|M(\mu_d)\| ~\leq ~ \sqrt{2} / \underline c_t.
\label{control:M_bound}
\end{equation}
We now focus our attention on the dynamics of the position error $e_p = p - p_r$ and the system velocity $v$.
Let $\tilde\mu = \mu - \mu_d$, where $\mu$ is the function defined by (\ref{control:v_dot}). In light of the choice for $\mu_d$, the derivatives of the position error and velocity can be written as
\begin{equation}
\dot e_p
~=~  v,
\qquad
\dot v
~=~ - k_p h(e_p) - k_v h(v) + \tilde\mu + \delta.
\label{control:v_dot:b}
\end{equation}
As previously mentioned, the velocity observer error  $\tilde v = v - \hat v$ is considered as a function of the apparent acceleration vector $r_2$. In fact, we define an error function associated to the apparent acceleration $r_2$ which is given by
\begin{equation}
\tilde r_2 =  k_1 \tilde v - ( I - \tilde R ) r_2.
\label{control:r2_tilde}
\end{equation}
Another important error function which we will focus on is the attitude error function $\tilde R$, or equivalently $\tilde Q = \left(\tilde\eta,\tilde q\right)$, which defines the relative orientation between the actual system attitude and the desired attitude. To prove the theorem, we will construct a Lyapunov function in terms of the error functions $\tilde q$, $\tilde r_2$, $v$ and $e_p$, in order to show that all of these states tend to zero. Since the dynamics of $\tilde q$ (or equivalently $\tilde\eta$), and $\tilde r_2$ are somewhat complicated, we will begin by first simplifying the expressions for their derivatives.
In order to analyze the dynamics of the attitude error, it is sufficient to study the derivative of the quaternion-scalar $\tilde\eta$. This is also desired since the derivative of the quaternion scalar can be less complicated than the derivative of the quaternion vector. As a starting point, the derivative of $\tilde\eta$ can be found from (\ref{control:Q_tilde_dot}) to be $\dot{\tilde\eta} = -\tilde q\tp\tilde\omega/2$ where $\tilde\omega = R_d\tp (\omega - \omega_d)$ and $\omega_d = M(\mu_d)\dot\mu_d$. To find a result for the desired angular velocity $\omega_d$ we first use the results (\ref{control:v_dot:b}), in addition to the derivative of the bounded function $h(\cdot)$, denoted as $\phi(\cdot)$ as defined in section (\ref{section:math_background}), to differentiate the virtual control law $\mu_d$ to obtain
$
\dot\mu_d =  - k_p \phi(e_p) v
- k_v \phi(v) \left( - k_p h(e_p) - k_v h(v)
+ \tilde\mu + \delta \right).
$
Simplifying this result, we obtain the following expression for the desired angular velocity
\begin{equation}
\omega_d = M(\mu_d) \left( f_{\mu_d}
- k_v \phi(v) \delta - k_v \phi(v) \tilde\mu \right).
\label{control:omega_d}
\end{equation}
Recall the control input $\omega$ uses the function $\psi$, given by (\ref{design:psi}). Using (\ref{control:r2_tilde}), the property (\ref{property:skew:so3}) and the fact $S(\tilde R r_2) \tilde R r_2 = 0$, $\psi$ can be rewritten as
\begin{equation}
\psi
= R_d \big(
 \gamma_1  S(r_1) \tilde R r_1
+ \gamma_2  S(r_2) \tilde R r_2
+ \gamma_2  S(\tilde r_2) \tilde R r_2
\big).
\label{control:psi:b}
\end{equation}
Finally, using the expression for the control input $\omega$, the error function $\tilde r_2$, in addition to (\ref{control:omega_d}), (\ref{control:psi:b})  and the fact $b_2 + u_t e_3 = R \delta$, we find the derivative
%
$
\dot{\tilde\eta }
=
-\frac{1}{2} \tilde q\tp R_d\tp\Big(
\gamma_1 R_d S(r_1)\tilde R r_1
+ \gamma_2 R_d S(r_2)\tilde R r_2
+ \gamma_2 R_d S(\tilde r_2) \tilde R r_2
 + k_v M(\mu_d) \phi(v)( I - \tilde R) \delta
+ k_v M(\mu_d) \phi(v) \tilde\mu
\Big).
$
%
To further simplify this result, we first recognize that in light of the definition of the rotation matrix from (\ref{q2r}) and the property $S(u)u = 0$, one can find
$
\tilde q\tp S(r_i) \tilde R r_i =
 2 \tilde q \tp S(r_i) \left( \tilde q \tilde q\tp - \tilde\eta S(\tilde q) \right) r_i
=
2 \tilde\eta \tilde q\tp S(r_i)^2 \tilde q
$.
Therefore, using the expression for the matrix $W$ defined by assumption \ref{assumption:lambda_min_w}, we obtain
\begin{align}
\dot{ \tilde\eta}
&=
\tilde\eta \tilde q\tp W \tilde q - \frac{\gamma_2}{2} \tilde q\tp S(\tilde r_2) \tilde R r_2
\nn \\ &
 -\frac{k_v}{2} \tilde q\tp R_d\tp M(\mu_d) \phi(v) \left( ( I - \tilde R ) \delta + \tilde\mu\right).
\label{proof:eta_e_dot:c}
\end{align}
Note that due to assumption \ref{assumption:lambda_min_w}, the matrix $W$ is positive-definite. We now shift our focus to study the dynamics of the error function $\tilde r_2$. In light of the expression for $\dot v$ from (\ref{control:v_dot:a}), the expression for $\dot{\hat v}$ from (\ref{design:v_hat_dot}), the attitude error dynamics from (\ref{control:Q_tilde_dot})-(\ref{control:omega_tilde}), the expressions (\ref{design:omega}), (\ref{control:omega_d}), and using the fact that $-k_1 \tilde v + r_2 - \hat R\tp b_2 = -\tilde r_2$,  we obtain
\begin{align}
\dot{\tilde r}_2
&=  - k_1 \tilde r_2 - ( I - \tilde R) \dot r_2
\nn \\ &
+ k_v R_d\tp S(b_2) M(\mu_d) \phi(v)  ( ( I - \tilde R)\delta + \tilde\mu ).
\label{control:r2_tilde_dot}
\end{align}
A commonality between the dynamic equations for $\dot{\tilde\eta}$ and $\dot{\tilde r}_2$, is that they both depend on the error functions $(I - \tilde R)$ and $\tilde\mu$. These two error functions can both be expressed in terms of the attitude error using the quaternion vector part $\tilde q$, which will be a useful characteristic later in the Lyapunov analysis. To describe this relationship we define two functions, $f_1(u_t,\tilde\eta,\tilde q),f_2(x,\tilde\eta,\tilde q) \in\mathbb{R}^{3\times 3}$ such that
\begin{equation}
\tilde\mu
~=~
  f_1(u_t,\tilde\eta,\tilde q) \tilde q,
\quad
( I - \tilde R )x
~= ~
f_2(x,\tilde\eta,\tilde q)\tilde q,
\end{equation}
where $x \in \mathbb{R}^3$. Using the definition of $\tilde\mu = \mu - \mu_d$, in addition to the expressions for $\mu$ and $\mu_d$ from (\ref{control:v_dot}) and (\ref{control:mu_d_defn}), respectively, one can find $f_1(u_t,\tilde\eta,\tilde q) = 2 u_t \left( \tilde\eta I - S(\tilde q) \right) S(R\tp e_3)$ and $f_2(x,\tilde\eta,\tilde q) = 2 ( S(\tilde q) - \tilde\eta I ) S(x)$. Based upon these definitions and the fact that $\|\tilde\eta I - S(\tilde q) \| = 1$, we find the following upper bounds for these two functions
\begin{equation}
\| f_1(u_t,\tilde\eta,\tilde q) \| \leq 2 \bar c_t,
~ \| f_2(x,\tilde\eta,\tilde q) \| \leq 2 \| x \|,
\label{proof:function_bounds}
\end{equation}
We now propose the following Lyapunov function candidate:
\begin{equation}
\V
=
\gamma k_p \left( \sqrt{1 + e_p\tp e_p} - 1 \right)
+ \frac{\gamma}{2} v\tp v
+ \frac{\gamma k_r}{2} \tilde r_2\tp\tilde r_2
+ \gamma_q \left( 1 - \tilde\eta^2\right),
\label{proof:V}
\end{equation}
\noi where $\gamma,\gamma_q,k_p$ and $k_r$ are positive constants. In light of (\ref{control:v_dot:b}), (\ref{proof:eta_e_dot:c}), (\ref{control:r2_tilde_dot}),  we have
\begin{align}
\dot\V
&=
- \gamma k_v v \tp h(v) + \gamma v\tp \delta - \gamma k_r k_1 \tilde r_2\tp\tilde r_2 - 2 \gamma_q \tilde\eta^2 \tilde q\tp W  \tilde q
\nn \\ &
+ \gamma k_v k_r \tilde r_2\tp R_d\tp S(b_2) M(\mu_d) \phi(v) \big( f_1(u_t,\tilde\eta,\tilde q)
\nn \\ &
+ f_2(\delta,\tilde\eta,\tilde q) \big) \tilde q
- \gamma k_r \tilde r_2 \tp f_2(\dot r_2,\tilde\eta,\tilde q) \tilde q
+ \gamma v\tp f_1(u_t,\tilde\eta,\tilde q) \tilde q
\nn \\ &
+ \gamma_q k_v \tilde\eta \tilde q\tp R_d\tp M(\mu_d) \phi(v) \left( f_1(u_t,\tilde\eta,\tilde q) + f_2(\delta,\tilde\eta,\tilde q) \right) \tilde q
\nn \\ &
+ \gamma_2 \gamma_q \tilde\eta \tilde q\tp S(\tilde r_2) \tilde R r_2.
\end{align}
Now, we wish to show that for an appropriate choice of the control gains, $\dot\V$ is guaranteed to be non-positive. However, this objective is a bit involved, and therefore requires we study the bound of several functions used in the expression of $\dot V$. We begin this analysis by defining the function scalar function $\sigma(t) := \sqrt{2 \V(t) }$.\\
 Based upon the definition of $\V$ from (\ref{proof:V}), the states $v$ and $\tilde r_2$ are bounded by $\sigma$ as follows
$
\| v(t) \| \leq \sigma(t) / \sqrt{\gamma}
,$
$
\| \tilde r_2(t) \| \leq \sigma(t) / \sqrt{\gamma k_r}
$.
Therefore, in light of assumption \ref{assumption:aerodynamic_forces}(\ref{assum:aero:b}), one can conclude that
\begin{equation}
\| \delta(v) \| \leq c_1 \sigma(t)^2 / \gamma.
\label{proof:delta_norm}
\end{equation}
\noi Due to the bounds of the functions $f_1(u_t,\tilde\eta,\tilde q)$ and $f_2(\delta,\tilde\eta,\tilde q)$ from (\ref{proof:function_bounds}), and the definition of $r_2$ from (\ref{control:v_dot:a}) we also find
\begin{align}
\| f_1(u_t,\tilde\eta,\tilde q) + f_2(\delta,\tilde\eta,\tilde q)\|
&\leq
2\left(\gamma \bar c_t + c_1 \sigma(t)^2\right)/ \gamma
\\
\| b_2 \|
&\leq
 \left( \gamma \bar c_t + c_1 \sigma(t)^2 \right) / \gamma.
\label{proof:b2_norm}
\end{align}
Given these bounds, we now apply Young's inequality to a number of the undesired terms in the expression for $\dot\V$:
\begin{equation}
\ba{l}
\gamma v\tp f_1(u_t,\tilde\eta,\tilde q) \tilde q
\\
\leq
\gamma \left( \Frac{1}{2} \left( \frac{\epsilon_1}{\sqrt{1+v\tp v}}\right) v\tp v + \frac{1}{2}\left(\frac{\sqrt{1 + v\tp v}}{\epsilon_1}\right) 4 \bar c_t^2 \tilde q\tp \tilde q
\right)
 \\
\leq
\Frac{\gamma \eps_1}{2} v\tp h(v) +
\Frac{2 \sqrt{\gamma} \bar c_t^2}{\eps_1}\sqrt{\gamma + \sigma(t)^2} \tilde q\tp \tilde q,
\ea
\end{equation}
\begin{equation}
\ba{l}
\gamma k_v k_r \tilde r_2\tp R_d\tp S(b_2) M(\mu_d) \phi(v) \left( f_1(u_t,\tilde\eta,\tilde q) + f_2(\delta,\tilde\eta,\tilde q) \right) \tilde q
\\
\leq
\Frac{\gamma k_v k_r \eps_2}{2} \tilde r_2\tp \tilde r_2
\\ \quad
+ \Frac{\gamma k_v k_r }{2\eps_2}
\left( \Frac{2}{\underline c_t^2}\right)
\left(\Frac{4\left(\gamma \bar c_t + c_1 \sigma(t)^2\right)^4}{\gamma^4}\right) \tilde q \tp \tilde q
\\
\leq
\Frac{\gamma k_v k_r \eps_2}{2}\tilde r_2\tp\tilde r_2
+ \Frac{4 k_v k_r}{\eps_2 \gamma^3 \underline c_t^2}
\left( \gamma \bar c_t + c_1 \sigma(t)^2 \right)^4 \tilde q\tp \tilde q,
\ea
\end{equation}
\noi where the norm of $M(\mu_d)$ is given by (\ref{control:M_bound}). To determine the bound of the term involving the time-derivative of $r_2$, we first derive the expression for $\dot r_2$ to be
\begin{equation}
\ba{l}
\dot r_2
=
- \dot u_t R\tp e_3 + u_t R\tp S(e_3) \omega + \dot\delta
 \\=
- \Frac{1}{u_t} \left(\mu_d - g e_3 \right)\tp
\Big( - k_p \phi(e_v) v - k_v \phi(v) f_1(u_t,\tilde\eta,\tilde q) \tilde q
 \\
+ k_v \phi(v) \left( k_p h(e_p) + k_v h(v) \right)
- k_v \phi(v) \delta\Big)R\tp e_3
\\
+ u_t R\tp S(e_3)
\bigg( M(\mu_d) \Big( - k_p \phi(e_v) v - k_v \phi(v) \tilde R \delta
\\
+ k_v \phi(v) \left( k_p h(e_p) + k_v h(v) \right)\Big) + \gamma_1 S(R_d r_1) b_1
 \\
 + \gamma_2 R_d S(\tilde r_2) \tilde R r_2
+ \gamma_2 R_d S(r_2) \tilde R r_2 \bigg) + \dot\delta.
\ea
\end{equation}
\noi Due to the bounds of the functions $h(\cdot)$, $\phi(\cdot)$, the (upper and lower) bounds of the thrust control input $u_t$, the bound of $\dot \delta$ from assumption \ref{assumption:aerodynamic_forces}(\ref{assum:aero:c}), the bound of $b_2$ from (\ref{proof:b2_norm}) (same as the bound of $r_2$), and the bound of $\delta$ from (\ref{proof:delta_norm}), we find that there exists five positive constants $d_{i}>0$, such that the norm of $\dot r_2$ is bounded by
$
\dot r_2 \leq
d_1 + d_2 \| v\| + d_3 \| v\|^2 + d_4 \| v\|^3 + d_5 \| v\|^4
$.
However, for the sake of simplicity, from this result we further conclude that there exists positive constants $c_3$ and $c_4$ such that
$
\dot r_2
\leq
c_3 + c_4 \sigma(t)^4.
$
As a result of this analysis, we again use Young's inequality to establish the following bounds:
\begin{equation}
\ba{l}
\gamma k_r \tilde r_2\tp f_2(\dot r_2,\tilde\eta,\tilde q) \tilde q
 \\
\qquad \leq
\Frac{\gamma k_r \epsilon_3}{2} \tilde r_2\tp \tilde r_2
+ \Frac{2 \gamma k_r}{\eps_3}
\left(
c_3 + c_4 \sigma(t)^4
\right)^2
\tilde q\tp \tilde q,
\ea
\end{equation}
\begin{equation}
\ba{l}
\gamma_2 \gamma_q \tilde\eta \tilde q\tp S(\tilde r_2) \tilde R r_2
\\  \leq
\Frac{\gamma_2^2 \gamma_q \eps_4}{2}\tilde r_2\tp\tilde r_2
+ \frac{\gamma_q}{2 \gamma^2\eps_4}
\left(
\bar c_t \gamma + c_1 \sigma(t)^2
\right)^2 \tilde\eta^2 \tilde q\tp \tilde q,
\ea
\end{equation}
\begin{equation}
\ba{l}
\gamma_q k_v \tilde\eta \tilde q\tp R_d\tp M(\mu_d) \phi(v) \left( f_1(u_t,\tilde\eta,\tilde q) + f_2(\delta,\tilde\eta,\tilde q) \right) \tilde q
\\\quad \leq
\gamma_q k_v \left( \Frac{\sqrt{2}}{\underline c_t}\right)
\left(
\Frac{2 \left( \gamma \bar c_t + c_1 \sigma(t)^2\right)}{\gamma}
\right) |\tilde\eta| \tilde q\tp \tilde q
 \\ \quad\leq
\Frac{2 \sqrt{2} \gamma_q k_v \left(
\gamma \bar c_t + c_1 \sigma(t)^2
\right)}{\gamma \underline c_t}
|\tilde\eta| \tilde q\tp \tilde q,
\ea
\end{equation}
\noi Recall from assumption \ref{assumption:lambda_min_w} that the norm of the matrix $W$ has a lower bound which is denoted as $c_w$. Therefore, in light of the lower bounds defined above, and assumption \ref{assumption:aerodynamic_forces}(\ref{assum:aero:a}) we find the expression $\dot \V$ is bounded by
\begin{equation}
\ba{l}
\dot\V(t)
 \leq
- \gamma v\tp h(v) \left( k_v - \eps_1/2 \right)
\\
- \gamma k_r \tilde r_2\tp \tilde r_2
\Big(
k_1 - \Frac{ \eps_2 k_v + \eps_3}{2}
- \Frac{\eps_4 \gamma_2^2 \gamma_q}{2 \gamma k_r}
\Big)
- \gamma_q \tilde\eta^2 \tilde q\tp \tilde q \Big( 2 c_w
\\
- \Frac{1}{\tilde\eta^2} \left(\Frac{\alpha_1(t)}{\eps_1}
+\Frac{\alpha_2(t)}{\eps_2}
+
\Frac{\alpha_3(t)}{\eps_3}\right)
-\Frac{\alpha_4(t)}{\epsilon_4}
 \\
- \Frac{2 \sqrt{\gamma} \bar c_t^2 \left( \gamma + \sigma(t)^2\right)^{1/2}}{\tilde\eta^2}
- \Frac{2 \sqrt{2} k_v \left( \gamma \bar c_t + c_1 \sigma(t)^2\right)}{\gamma \underline c_t |\tilde\eta|}
\Big),
\ea
\label{proof:V_dot}
\end{equation}
\begin{align}
\alpha_1(t)
&=
2 \sqrt{\gamma} \bar c_t^2 \sqrt{\gamma + \sigma(t)^2} / \gamma_q,
\\
\alpha_2(t)
&=
4 k_v k_r \left(
\gamma\bar c_t + c_1\sigma(t)^2\right)^4 / (
\gamma^3 \underline c_t^2 \gamma_q),
\\
\alpha_3(t)
&=
2 \gamma k_r
\left(
c_3 + c_4 \sigma(t)^4
\right)^2 /  \gamma_q,
\\
\alpha_4(t)
&=
 \left(
\bar c_t \gamma + c_1 \sigma(t)^2
\right)^2 / (2 \gamma^2).
\end{align}
Now, let us define a lower bound for $|\tilde\eta|$, which based upon some appropriate choices of gains, ensures $\dot\V \leq 0$ for all $t \geq t_0$.
Note that when $\tilde\eta(t) = 0$ we cannot guarantee stability using (\ref{proof:V_dot}) since in this case $\dot\V$ could potentially be positive. To show that $\tilde\eta(t)$ is never zero, we first introduce the positive constant $\rho$ which is the desired lower bound for $|\tilde\eta(t)|$. Therefore, $\rho$ must be chosen to satisfy
$
0 < \rho < |\tilde\eta(t_0)|.
$
Subsequently, based upon the definition of the Lyapunov function candidate (\ref{proof:V}), we choose
$
\gamma = \bar\gamma (
k_p (\sqrt{ 1 + \|e_p(t_0)\|^2}-1) + \|v(t_0)\|^2/2
+ \|\tilde r_2(t_0)\|^2/2 + \xi
)^{-1},
$,
%
where the parameter $\xi$ is chosen to be positive, and $\bar\gamma$ is chosen to satisfy
%
$
0 < \bar\gamma <  \gamma_q \left( \tilde\eta(t_0)^2 - \rho^2 \right),
$
%
where $\gamma_q$ is chosen to be positive. Recall $k_p>0$ and $k_v>0$ are chosen arbitrarily provided that $k_p + k_v < g$. The remaining gains and parameters are chosen to ensure that all terms in (\ref{proof:V_dot}) are guaranteed to be negative at the initial time $t_0$. The gains and parameters are chosen as follows: Choose $\epsilon_1$ such that $0 < \epsilon_1 < 2 k_v$. Recall that the minimum eigenvalue of $W$, denoted by $c_w>0$, can be increased using the gains $\gamma_1$ and $\gamma_2$. Therefore, there exists constants $\bar\gamma_1$,$\bar\gamma_2$, and $\bar\eps_i$, $i = 2,3,4$, such that for all $\gamma_1 > \bar\gamma_1$, $\gamma_2 > \bar\gamma_2$, and $\eps_i > \bar\eps_i$, the following inequality is satisfied
\begin{equation}
\ba{l}
2 c_w >
 \Frac{1}{\rho^2} \left(\Frac{\alpha_1(t_0)}{\eps_1}
+\Frac{\alpha_2(t_0)}{\eps_2}
+
\Frac{\alpha_3(t_0)}{\eps_3}\right)
+\Frac{\alpha_4(t_0)}{\epsilon_4}
\\
+\Frac{2 \sqrt{\gamma} \bar c_t^2 \left( \gamma + \sigma(t_0)^2\right)^{1/2}}{\rho^2}
+ \Frac{2 \sqrt{2} k_v \left( \gamma \bar c_t + c_1 \sigma(t_0)^2\right)}{\gamma \underline c_t \rho}.
\ea
\label{c_w_inequality}
\end{equation}
Finally, choosing
$k_1 > \kappa_1(\epsilon_2,\epsilon_3,\epsilon_4,\gamma)
:= (\eps_2 k_v +\eps_3)/2
+(\eps_4 \gamma_2^2 \gamma_q) / (2 \gamma k_r)
$
we conclude that $\dot\V(t_0) \leq 0$ at the initial time $t_0$. We now need to show that this is true for all time. Since the functions $\alpha_1(t)$ through $\alpha_4(t)$ are non-increasing if $\dot V \leq 0$, then a sufficient condition for $\dot\V(t) \leq 0$ is $|\tilde\eta(t) | \geq \rho$. We will now show that indeed $\rho \leq |\tilde\eta(t)|$ for all $t > t_0$. Suppose that there exists a time $t_1$ such that for all $t_0 \leq t < t_1$, $|\tilde\eta(t)| \geq \rho$ and $|\tilde\eta(t_1)| < \rho$ when $t = t_1$. At the time $t_1$ from (\ref{proof:V}), it is clear that
$
\V(t_1)
~\geq~
 \gamma_q \left( 1 - \tilde\eta(t_1)^2 \right)
~>~ \gamma_q \left( 1 - \rho^2 \right).
$
However, due to the choice of $\gamma$ and $\bar\gamma$ the value of the Lyapunov function candidate at the initial time $t_0$ must satisfy
$
\V(t_0)
 <
\bar\gamma + \gamma_q \left( 1 - \tilde\eta(t_0)^2\right)
<
\gamma_q \left( 1 - \rho^2 \right)
$
\noi and therefore $\V(t_1) > \V(t_0)$. This is a contradiction since $\dot\V(t) \leq 0$ for all $t_0 \leq t < t_1$, and the functions $\V(t)$, $\alpha_i(t)$ and $\sigma(t)$ are non-increasing in the interval $t_0 \leq t < t_1$. Therefore, we conclude that $|\tilde\eta(t)| \geq \rho$ and $\dot\V(t) \leq 0$ for all $t > t_0$, and the states $v$ and $\tilde r_2$ are bounded. Therefore, $\dot{\tilde r}_2$, $\dot v$, $\dot{\tilde{\eta}}$, and  $\ddot\V$ are bounded. Invoking Barbalat's Lemma, one can conclude that $\lim_{t\to\infty}\left(v(t),\tilde r_2(t),\tilde q(t)\right) = 0$. Furthermore, since $\lim_{t\to\infty} \dot v(t) = 0$, and $\lim_{t\to\infty} \delta(t) = 0$, it follows from the expression of the velocity dynamics $\dot v = - k_p h(e_p) - k_v h(v) - \delta = 0 $, that $\lim_{t\to\infty} e_p(t) = 0$, which ends the proof. $\hfill\blacksquare$
%
%
%
%
%
%
\section{Simulations}
\begin{figure}
\centering
\includegraphics[width = 0.7\columnwidth]{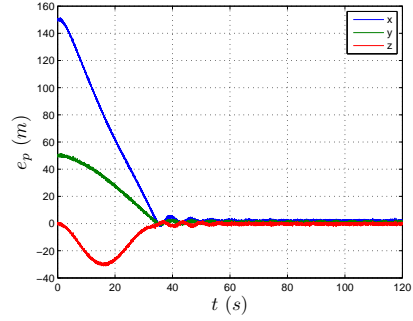}
\caption{Position Error $~e_p~ (m)$}
\label{figure:control:position_error}
\end{figure}

Simulation results have been provided for the system defined by (\ref{control:p_dot})-(\ref{control:Q_dot}) with the proposed control laws (\ref{design:omega})-(\ref{design:mu_d}), which are shown in Figure (\ref{figure:control:position_error}). To demonstrate the robustness of the proposed control strategy, we have included in our simulations the presence of wind, sensor noise and gyro-bias. A much more aggressive aerodynamic model is adopted during the simulations, which violates the assumptions in order to demonstrate the robustness of the proposed controller. This aerodynamic model considers that the aerodynamic drag of the system is a function of the system attitude, and also considers that the system is operating in the presence of uniform wind. The following model was used for the aerodynamic disturbance $\delta$:
$
\delta = - \frac{1}{m_b} \| v - v_w \| R\tp C_d R (v - v_w)
$
\noi where $v_w \in \mathbb{R}^3$ is the inertial referenced wind velocity vector, $m_b$ is the system mass and $C_d = C_d\tp > 0$ is a constant positive definite matrix that represents body-referenced aerodynamic drag coefficients that are dependent on the system geometry. Note that the expression for $\delta$ depends on the mass $m_b$ of the system due to the definition of the velocity dynamics from (\ref{control:v_dot}). A constant wind velocity vector was specified as $v_w = [10,5,0] m/s$. For this simulation the value $C_d = \mbox{diag}\left[0.1,0.1,0.05\right] kg/m$ was chosen and the system mass was specified as $m_b = 5kg$.  Gaussian noise was added to the magnetometer, accelerometer, gyroscope, linear velocity and position sensors with standard deviation values equal to $0.01G$, $0.1m/s^2$, $0.1 deg/s$, $0.5m/s$ and $0.5m$, respectively. A bias was also added to the gyro measurements which was chosen as $[0.1,0.05,-0.2] deg/s$. The system gains were chosen as follows: $k_p = 5$, $k_v = 0.1$, $\gamma_1 = 0.1 = \gamma_2 = 0.05$ and  $k_1 = 5$. The following initial conditions were used: $p(t_0) = \left[150,50,0\right]m$, $v(t_0) = \left[0,0,0\right]m/s$, $\hat v(t_0) = \left[0,0,0\right]$, $Q(t_0) = \left[ 1,0,0,0\right]$, and the desired position was chosen as $p_r = \left[0,0,0\right]$. The inertial vector for the magnetometer measurement was chosen as $r_1 = [0.18,0,0.54]G$, which is consistent with the magnitude and direction of the Earths magnetic field in Southern Ontario, Canada.
The simulation results show that the system position converges to the desired value, except for a small error attributed to the sensor noise, gyroscope bias and the effect of the aerodynamic disturbance caused by wind.

\begin{figure}
\centering
\includegraphics[width =0.7\columnwidth]{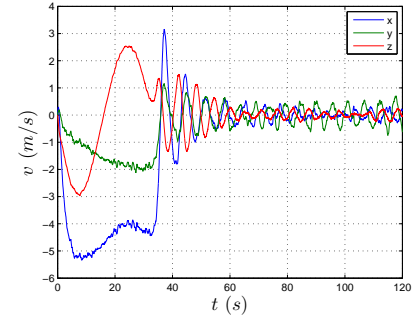}
\caption{Velocity  $~v (m/s)$}
\label{figure:control:velocity}
\end{figure}

\begin{figure}
\centering
\includegraphics[width = 0.7\columnwidth]{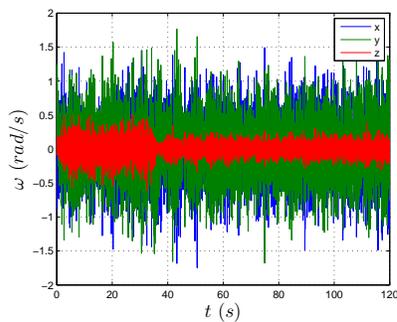}
\caption{Control Effort (Angular Velocity) $~\omega~$ }
\end{figure}

\begin{figure}
\centering
\includegraphics[width = 0.7\columnwidth]{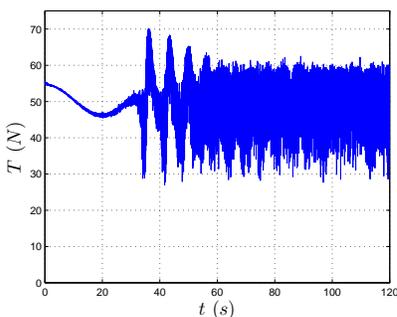}
\caption{Control Effort (Thrust) $~T = u_t \cdot m_b~ (N)$}
\label{figure:control:thrust}
\end{figure}
\section{Conclusion}
A new position controller for VTOL-UAVs that does not require direct measurement of the system's attitude, nor does it require the use of an attitude observer has been proposed. The accelerometer and magnetometer signals are explicitly used in the control law to capture the necessary information about the system's orientation (without, explicitly, reconstructing or estimating the orientation). Furthermore, the usual simplifying assumption restricting the accelerometer measurement to the gravity direction is not required anymore. In fact, in this work, the accelerometer is used to measure the system's apparent acceleration, which makes the proposed control efficient when the system is subjected to significant linear accelerations. We have shown that, through an appropriate choice of the control gains, the system position is guaranteed to be bounded and to converge to the desired position for almost all initial conditions. Simulation results have also been performed which demonstrate the effectiveness of the controller,
despite the choice of relatively low control gains, and large initial conditions.
\begin{ack}                               
This work was supported by the National Sciences and Engineering Research Council of Canada (NSERC)  
\end{ack}

\bibliographystyle{apalike}        
\bibliography{Roberts_Tayebi_2011}  



\end{document}